\input amstex
\documentstyle{amsppt}
\magnification=\magstep1
\pagewidth{6.2truein}
\pageheight{8.5truein}
\hcorrection{5mm}
\vcorrection{12mm}
\baselineskip=14pt
\parindent=12mm

\def\BC{\Bbb C}
\def\BS{\Bbb S}

\def\BB{\Bbb B}
\def\BD{\Bbb D}
\def\hat{\widehat}
\def\tilde{\widetilde}
\def\epsilon{\varepsilon}
\def\grad{\overline\nabla_X f}

\topmatter
\title  EULER OBSTRUCTION AND DEFECTS OF FUNCTIONS ON SINGULAR VARIETIES
\endtitle
\author  J.-P. Brasselet, D. Massey, A. J. Parameswaran and J. Seade  
\endauthor
\address Jean-Paul Brasselet: Institut de Math\'ematiques de Luminy,
UPR 9016 CNRS, Campus
de Luminy - Case 930, 13288 Marseille Cedex 9, France
\endaddress
\email jpb$\@$iml.univ-mrs.fr\endemail
\address David Massey: Department of Mathematics, Northeastern  
University,
567 Lake Hall,
Boston, MA   02115  USA, http://www.massey.math.neu.edu/
\endaddress
\email dmassey$\@$neu.edu \endemail
\address A. J. Parameswaran: Tata Institute of Fundamental Research,
Homi Bhaba Road, Colaba, Mumbai, India
\endaddress
\email param$\@$math.tifr.res.in\endemail
\address Jos\'e Seade:
Instituto de Matem\'aticas, Unidad Cuernavaca,
Universidad Nacional Aut\'onoma de M\'exico,
Apartado Postal 273-3, C.P. 62210, Cuernavaca, Morelos,
M\'exico
\endaddress
\email jseade$\@$matem.unam.mx
\endemail

\abstract Several authors have proved Lefschetz type formulae for the  
local Euler obstruction. In particular, a result of this type is proved in  \cite  
{BLS}. The formula proved in that paper turns out to be
equivalent to saying that the local Euler obstruction, as a  
constructible function, satisfies the local Euler condition (in bivariant theory) with respect  
to general linear forms. The purpose of this work is to understand what prevents the  
local Euler obstruction of satisfying the local Euler condition  with
respect to functions which are singular at the considered point. This is
measured by an invariant (or ``defect'') of such functions that we define below.
We give an interpretation of this defect in terms of vanishing cycles,  
which allows us to calculate it algebraically,
using \cite{Ma2, Ma3}. When the function has an isolated singularity, our  
invariant can be defined geometrically, via obstruction theory.
We notice  that this invariant  unifies the usual concepts of {\it the Milnor number} of a  
function and of the {\it  local Euler obstruction} of an analytic set.

\endabstract

\thanks{\hskip-20pt Partially supported  by the Cooperation Programs
France-M\'exico CNRS/CONACYT, France-India CNRS/NBHM, and by the
CONACYT Grant G36357-E.}
\endthanks



\endtopmatter
\document
\NoRunningHeads
\head
{0. INTRODUCTION.}
\endhead

This work is a natural continuation of \cite {BLS},  where it is proved  
a
Lefschetz-type formula for the local Euler obstruction.
More precisely, let $ (X,0) $  be an equidimensional complex analytic
singularity germ in an open set $ U \subset \BC^N$. We endow
$ (X,0) $ with a Whitney stratification  $\{ V_i\}$ and
consider a general  complex linear form $l : \BC^N \to \BC$. The  
formula of
\cite {BLS} says that the local Euler obstruction of $X$ at $0$  
satisfies:
$$Eu_X(0)=\sum_i  \chi(V_i\cap \Bbb B_\varepsilon\cap
l^{-1}(t_0)) \cdot Eu_X(V_i),$$
where $\Bbb B_\varepsilon$ is a small ball around $0$ in $\BC^N$, $t_0  
\in
\BC
\setminus \{ 0 \}$ is sufficiently near $\{ 0 \}$ and $Eu_X(V_i)$ is the
Euler
obstruction of $X$ at any point of the stratum $V_i$.

\medskip
As we explain  below, this formula is actually equivalent to  saying  
that
the local
Euler
obstruction, as a constructible function on $X$, satisfies the local  
Euler
condition  with respect to general linear forms.
The purpose of this work is to understand what prevents the local Euler
obstruction of satisfying the local Euler condition  with
respect to functions with a singularity at $0 \in X$.
\medskip

We consider first an analytic
function  $f: X \to \BC$ with an isolated singularity at $0$.
We introduce (in section 2) an invariant $Eu_{f,X}(0)$, called
{\it the local Euler obstruction} of $f$ at $0$.
To define this invariant, let us consider the Nash bundle $\tilde T$
over the Nash transform  $\tilde X$ of $X$.
Roughly speaking,
$Eu_{f,X}(0)$ is the obstruction for constructing a section of
$\tilde T$ that extends a lifting of $\grad$,
the complex conjugate of the gradient vector field of $f$ on $X$ (see
section 2 for a precise definition). In section 3, we prove (Theorem  
3.1):
$$Eu_X(0)= \Big(\sum_i \chi(V_i\cap \Bbb B_\varepsilon\cap
f^{-1}(t_0)) \cdot Eu_X(V_i)\Big) \, +\, Eu_{f,X}(0)\,  .$$
\medskip

In sections 3 and 4, we give a
proof of this formula, similar to the proof in \cite{BLS}.
  The main difference is that the invariant $Eu_{f,X}(0)$ does not
vanish in general and its contribution to the local Euler obstruction  
$Eu_X(0)$
has to be taken into account.
The more natural way to prove our formula is to work directly on $X$,  
which
is somehow
simpler than working on the Nash transform $\widetilde X$,
as in \cite {BLS}. It seems to us that this approach gives interesting
insigths even for the case of general linear forms.
\medskip

In section 5 we look at the situation where the function $f$ may have a  
non-isolated
singularity at $0$. In this case our geometric definition of the  
invariant
$Eu_{f,X}(0)$ does not make sense. However, one can define the following
{\it defect}:
$$D_{f, X}(0):=Eu_X(0) - \Big(\sum_i\chi\big(V_i\cap B_\epsilon
\cap f^{-1}(t_0)\big)\cdot Eu_X(V_i)\Big)$$
This is well-defined even when
the singularity is non-isolated, and
we use a formal, derived category argument to prove that the defect
$D_{f, X}(0)$
has a nice interpretation in terms of vanishing cycles  (Theorem 5.3).
The proof that we give is very short,
but uses a number of basic results, including the index theorem of
Brylinski, Dubson, and Kashiwara \cite{BDK}.
 From the defect point of view, what Theorem 3.1 tells us is that,
in the case of isolated singularities,  $D_{f, X}(0)= Eu_{f, X}(0)$.
\medskip

This vanishing cycle interpretation that we give for $D_{f, X}(0)$  
allows
us to calculate it algebraically. 
In the case of isolated singularities,
we may use another index theorem -- proved independently by Ginsburg
\cite{Gi}, L\^e \cite{Le1}, and Sabbah \cite{Sa2} -- to conclude
that, up to sign, $D_{f, X}(0)$ is equal to the intersection number
of the conormal variety of $X$  with the image of the differential
$d\tilde f$ at the point $(0, d_0\tilde f)$ (Corollary 5.4).
\medskip

Thus, when $X$ is
$\Bbb C^d$ and $f$ has an isolated singularity at $0$, we obtain
  that  $D_{f, X}(0)=(-1)^d\mu$. That is, $D_{f, X}(0)$
coincides (up to sign) with the usual Milnor
number of $f$ (compare Remark 3.4). We can also define the invariant
  $Eu_{f,X}(0)$ when $f$ is the function distance to $0$ on $X$; in this  
case
it  coincides with the usual local Euler number
$Eu_X(0)$. Hence $Eu_{f,X}(0)$ unifies these two important notions. The  
definition of the defect
$D_{f, X}(0)$ also makes sense when $f$ is a constant function, and in
this case the algebraic interpretation says that $D_{f, X}(0) =  
Eu_X(0)$.

\medskip

Moreover, 
even when $f$ has non-isolated singularities, in section 6 
we obtain algebraic formulas for $D_{f, X}(0)$ by using
the work of Massey. When $X$ is affine space, we use the results
of \cite{Ma2} on L\^e numbers. When $X$ is arbitrary, we use
the recent results of \cite{Ma3} on L\^e-Vogel numbers.
\medskip

It would be interesting
to understand
the relationship of this invariant with the {\it Milnor number} for  
functions on
singular spaces, introduced by V. Goryunov in \cite {Go} for functions  
on
curves (see
also
\cite {MS}). This latter invariant  has an obvious
generalization to functions on isolated complete intersection germs,  
either
as the number
of critical points (on a Milnor fibre) of a generic perturbation of the
given function,
or equivalently, as the GSV-index of the gradient vector field of the  
given
function \cite {GSV}.

\bigskip

\head
{1. EULER OBSTRUCTION AND THE EULER CONDITION.}
\endhead

Let us recall the result of \cite{BLS}.
Let $ (X,0) $  be an equidimensional complex analytic singularity germ  
of
dimension $d$
in an open set $ U \subset \BC^N$. Let $G(d,N)$ denote the Grassmanian  
of
complex $d$-planes in $\BC^N$. On the regular part $X_{reg}$ of $X$,  
there
is
a map $\sigma : X_{ reg} \to U\times G(d,N)$, defined by
$\sigma(x) = (x,T_x(X_{reg}))$. The Nash transformation $\tilde X$ of  
$X$
is the closure of Im$(\sigma)$ in $
U\times G(d,N)$. It is a complex analytic space endowed with an
analytic projection $\nu : \tilde X \to X$, which is a biholomorphism
on $\nu^{-1}(X_{reg})$. Let us denote by  $U(d,N)$ the
tautological bundle over $G(d,N)$, and by ${\Bbb U}$ the
corresponding locally trivial extension bundle over $ U \times G(d,N)$.
The Nash bundle $\tilde T$ on $\tilde X$ is the
restriction of  ${\Bbb U}$ to $\tilde X$. \par

Consider a complex analytic stratification $\{ V_i\}$  of
$U\subset \BC^N$, adapted to $X$. We
assume that $\{0\}$ is a stratum and
the stratification satisfies the Whitney conditions \cite{Wh}.
We choose a representative of the germ of $(X,0)$
sufficiently small, so that $\{0\}$ is in the closure of all the strata.

A {\sl stratified vector
field} $v$ on $X$ is a continuous section of the restriction
$TU\vert_X$ of $TU$ to $X$,
such that if $x \in V_i\cap X$, then $v(x) \in T_x(V_i)$.
Such a vector field is {\sl radial  } at $a \in X$
if for each sufficiently small ball  $\Bbb B_\varepsilon (a) \subset U $
around $a$, with $\varepsilon >0$, and for each boundary point
$x \in \Bbb S_\varepsilon (a) =\partial \Bbb B_\varepsilon (a)$, the
vector
$v(x)$ is pointing outwards the ball $\Bbb B_\varepsilon(a)$. Such a
vector
field has necessarily an isolated zero at $a$.

One has the following Lemma (\cite {BS}, Proposition 9.1):

\proclaim{ LEMMA 1.1} Let $v$  be a stratified non-zero vector field on
$A\subset X$.  Then $v$  can be lifted to a section $\tilde v$ of  
$\tilde
T$ along $\nu^{-1} (A)$.
\endproclaim

\proclaim
{PROPOSITION-DEFINITION 1.2} $($\cite {BS}, Proposition 10.1; see \cite
{MP} for the original definition$)$. Let $a$ be a point in $X$ and $v$  
a radial vector
field defined
on $X \cap \Bbb B_\varepsilon (a)$.
Let $\tilde v$ be the lifting of $v$ to
$\nu^{-1} (X \cap \Bbb S_\varepsilon (a))$ as a section of $\tilde
T$.  The
local Euler obstruction $Eu_X(a)$ is   the obstruction to the extension
of
$\tilde v$ as a nowhere zero section of $\tilde T$ inside $\nu^{-1}  
(X\cap
\Bbb
B_\varepsilon (a))$.
\endproclaim

More precisely, let ${\Cal O} (\tilde v) \in
H^{2d}\big(\nu^{-1}(X \cap
\Bbb B_\varepsilon(a)), \nu^{-1}(X \cap \Bbb S_\varepsilon(a))\big)\,$  
be
the
obstruction cocycle for  extending  $\tilde v$ as a nowhere zero section
of
$\tilde T$ inside $\nu^{-1} (X\cap \Bbb B_\varepsilon(a))$.
The local Euler obstruction $Eu_X(a)$ is the evaluation of ${\Cal O}
(\tilde v)$  on the fundamental class of the pair $\big(\nu^{-1}(X \cap
\Bbb B_\varepsilon(a)), \nu^{-1}(X \cap \Bbb S_\varepsilon(a))\big)$.  
The
Euler obstruction is an integer. It is constant on each stratum of the
Whitney stratification.
\smallskip

We consider a complex
linear form $l : U \to \BC$, so that $0 \in l^{-1}(0)$, and we assume  
that
the kernel of $l$, $ker(l)$, is transversal to every
  generalized tangent space $T$ at $0$, {\it i.e.} to every limit
of tangent spaces $T_{x_n} (V_i)$, for every $V_i$ and  every sequence
$x_n\in V_i$ converging to $0$. Such a
linear form is said to be {\it general} (with respect to $X$). From now
on,
we will denote by $\BB_\varepsilon$ the ball $\BB_\varepsilon(0)$.

\medskip

\proclaim{THEOREM 1.3} {\rm (\cite {BLS} Theorem 3.1)} Let $(X,0)$  and  
the
Whitney
stratification
$\{ V_i\}$   be as above, and let $l: U \to \BC$ be a general linear  
form.
Then:
$$Eu_X(0)=\sum_i \chi(V_i\cap \Bbb B_\varepsilon\cap
l^{-1}(t_0)) \cdot Eu_X(V_i),$$
where $\varepsilon$ is sufficiently small, $t_0  \in \BC \setminus  
\{0\}$
is
close to $0$  and
$Eu_X(V_i)$ is the Euler obstruction of $X$ at any point of the stratum
$V_i$.
\endproclaim

This theorem can be stated through the framework of bivariant theory
\cite{FM}. We follow the formulation given in \cite {Br}.
For this  we recall  that a
function $\alpha : X \to
\Bbb Z$ is {\it constructible} if for each $n \in \Bbb Z$, the set
$\alpha^{-1}(n)$ is   constructible, {\it i.e.} it is obtained as a  
finite
number of
unions, intersections and differences of analytic subsets of $X$. Given  
a
constructible function
$\alpha : X \to \Bbb Z$  and
$A\subset X$,  the weighted Euler-Poincar\'e characteristic
$\chi(A;\alpha)$ is defined
by:
$$\chi (A;\alpha) : = \sum_{n\in \Bbb Z} n \cdot\chi (A\cap
\alpha^{-1}(n)).$$

If $\{ V_i\}$ is a stratification of $X$ such that the value of $\alpha$
is
constant on each stratum, we denote by $\alpha(V_i)$ its value in any point
of
$V_i$. We have, equivalently,
$$\chi (A;\alpha) = \sum_i \chi (A\cap V_i)\cdot\alpha(V_i)\,.$$

Given an analytic map $f : X \to {\Bbb C}$ and a constructible function
$\alpha : X \to \Bbb Z$,
we say that $\alpha$ {\it satisfies the local Euler condition
with respect to} $f$ if  for each $x \in X$, we have
$$\alpha(x) = \chi( \Bbb B_\varepsilon(x)\cap f^{-1}(t);\alpha)\,,$$
where $\Bbb B_\varepsilon(x)$
is  a small ball in $U \subset \BC^N$, $t \in \Bbb D(f(x)) \setminus
\{f(x) \}$ and  $\Bbb D(f(x))$ is a small disc in ${\Bbb C}$ centered  
at $f(x)$
(see \cite{Br}, \cite{Sa1}).

  It is clear  that Theorem 1.3 can be restated as:

\proclaim{THEOREM 1.4} The local Euler obstruction, as a constructible
function, satisfies the local Euler condition with respect to   general
linear forms.
\endproclaim

\vskip.2cm
\noindent {\bf REMARK 1.5}.
Recently, J. Sch\"urmann \cite {Scu} proved
that,  with the same hypotheses as in Theorem 1.3, if  a constructible
function
   $\alpha$  on $X$ is a linear combination of the functions
$Eu_{{\overline V_j}}(*)$,  then one has
$$\alpha(0)=\sum_i \chi(V_i\cap \Bbb B_\varepsilon\cap l^{-1}(t_0))  
\cdot
\alpha(V_i)\,.$$
  In other words, $\alpha$ satisfies the local Euler condition with
respect to
general linear forms. The proof in \cite{Scu} uses the index theorem for
the
vanishing cycle functor. We notice that this result  can also be
proved using 1.3, applying it to the closure of each stratum $\overline
V_j$.

\bigskip

\head
{2. LOCAL EULER OBSTRUCTION OF A FUNCTION}
\endhead

Let us recall some well known concepts about singularity theory,
which originate in the work of R. Thom.

Let $(X,0)$ be a complex analytic germ as above, contained in an open
subset $U$
of  $\BC^N$ and  endowed with a complex analytic Whitney stratification
$\{ V_i\}$. We assume further that $\{0\}$ is a stratum and that every
stratum
contains $0$ in its closure. For every point $x\in X$, we will denote by
$V_i(x)$ the stratum containing $x$. Let $f:X \to
\BC$ be a holomorphic function, which is the restriction of a  
holomorphic
function
$\tilde f : U \to \BC$. We recall \cite {GM} that a {\it critical point}
of
$f$ is a point $x \in X$ such that
$d \tilde f(x)({T_x (V_i(x))}) = 0$.  We say, following \cite  {Le1},
\cite
{GM}, that $f$  has an isolated singularity at $0 \in X$ relative to the
given
Whitney stratification, if
$f$ has no critical points in a punctured neighbourhood of $0$ in $X$.

Let us denote by $\overline\nabla \tilde f(x)$
the gradient vector field of $ \tilde f$ at a point $x \in U$, defined  
by
$\overline \nabla \tilde f(x) := (\frac {\overline{\partial\tilde
f}}{{\partial
x_1}},..., \frac {\overline{\partial\tilde  f}}{\partial x_N}) \,,$  
where
the bar
denotes  complex conjugation. From now on we assume that $f$ has an
isolated singularity
at $0\in X$. This implies that the kernel
${ker} (d \tilde f)$ is transverse to $T_x(V_i(x))$ in any point $x\in
X\setminus
\{ 0\}$. Therefore at each point $x\in
X\setminus
\{ 0\}$, we have:
$$Angle \langle \overline\nabla\tilde f(x), T_x(V_i(x))\rangle < \pi/2\,
,$$
so the projection of $\overline\nabla\tilde f(x)$ on $T_x(V_i(x))$,
denoted by
$\hat \zeta_i(x)$, is not zero.

Let $V_j$ be a stratum such that $V_i \subset \overline V_j$,
and let $\pi : \Cal U_i \to V_i$ be a tubular neighbourhood of
$V_i$ in $U$. Following the construction of M. H. Schwartz in (\cite  
{Sc},
\S 2),
we see that the Whitney
condition $(a)$  implies that at each point $y\in V_j\cap \Cal U_i$, the
angle of
$\hat
\zeta_j(y)$ and of the parallel extension of $\hat \zeta_i(\pi(y))$  is
small.  This property implies that these two vector fields are homotopic
on the boundary of $\Cal U_i$. Therefore, we can glue together
the vector fields $\hat \zeta_i$ to obtain a   stratified
vector field on $X$, denoted by $\grad$. This vector field is
homotopic
  to $\overline \nabla \tilde f\vert_X$ and  one has $\grad
\ne 0 $ unless $x = 0$.

Intuitively, what we are doing in the construction of $\grad$ is to  
take,
for each stratum $V_i$ of $X$, the gradient vector field of the  
restriction
of $f$ to $V_i$, and then gluing all  these vector fields together. For
instance,
we can  make the same construction taking  $f$ to be the
function distance to $0$, then $\grad$ is a radial vector field as in  
1.2.

\medskip

\noindent{\bf DEFINITION 2.1.}
Let $\nu : \tilde X \to X$ be the Nash transform of $X$.
Let $\tilde\zeta$ be the lifting of $\grad$ as a section of the Nash
bundle
$\tilde T$ over $\tilde X$ without singularity over $\nu^{-1}(X\cap
\Bbb S_\varepsilon)$. Let ${\Cal O} (\tilde\zeta)
\in H^{2d}\big(\nu^{-1}(X \cap
\Bbb B_\varepsilon), \nu^{-1}(X \cap \Bbb S_\varepsilon)\big)$ be the
obstruction cocycle
  for the extension of  $\tilde \zeta$
as a nowhere zero section of $\tilde T$ inside $\nu^{-1} (X\cap \Bbb
B_\varepsilon)$.
We define {\it the local Euler obstruction of $f$ on} $X$ at
$0$, denoted $ Eu_{f,X}(0)$, to be the evaluation of ${\Cal O}
(\tilde\zeta)$
on the fundamental class of the pair $(\nu^{-1}(X \cap
\Bbb B_\varepsilon), \nu^{-1}(X \cap \Bbb S_\varepsilon))$.

\medskip
Notice that $Eu_{f,X}(0)$ is an integer.
We remark that a reason for considering the conjugate gradient vector
field
$\grad$, and not the usual gradient vector field $\nabla_X f$, is
given by the following Lemma:

\proclaim{LEMMA 2.2}
The vector field
$\grad$ is the lifting, up to homotopy, of a vector
field on $\Bbb C$, via $d\tilde f$.
\endproclaim

\noindent
{\bf Proof.} The gradient vector field satisfies
$$ d\tilde f(\overline \nabla \tilde  f(x)) =
||\overline \nabla \tilde  f(x)||^2  \, \in \Bbb R  \setminus \{0\}
\qquad \text{ for } x \in X  \setminus \{0\}\,,$$
this means that it is the lifting, up to scaling, of a constant vector
field on a small  disk $\BD_\eta\subset \BC$.
\qed

  \medskip
It is easy to see that if $0$ is a smooth point of $X$ and also
a regular point of $f$, then $Eu_{f,X}(0)=0$. In the proposition 2.4
below, we prove that this is the case in a more general situation.

  \medskip
  \noindent
{\bf DEFINITION 2.3.} Let $(X,0) \subset (U,0)$ be a germ of analytic  
set
in
$\BC^N$ equipped with a Whitney stratification and let $f: (X,0)
\to (\BC,0)$ be an analytic map, restriction of a regular
holomorphic function $\tilde f:
(U,0) \to (\BC,0)$. We say that $0$ {\it is a general point of} $\,f$ if
the hyperplane $\,ker \,d\tilde f(0)\,$ is transversal in $\BC^N$ to  
every
generalized tangent space at $0$, {\it i.e.} to every limit
of tangent spaces $T_{x_n} (V_i)$, for every $V_i$ and  every sequence
$x_n\in V_i$ converging to $0$.

\medskip
We notice that for every $f$ as above, the general points of $f$ form a
non-empty open set on each
(open) stratum of $X$, essentially by Sard's theorem.
We also remark that Definition 2.3 provides a coordinate free way of  
looking at
general linear forms.  In fact the previous definition is equivalent to
saying
that, with an appropriate  local change of coordinates $\tilde f$ is a
linear
form in $U$, and it is general with respect to $X$ (compare Lemma 1.3 in
\cite{BLS}).

\medskip

\proclaim{PROPOSITION 2.4 }  Let $0$ be a general point of $f: (X,0) \to
(\BC,0)$. Then
$$Eu_{f,X}(0) = 0 \,.$$
\endproclaim

\noindent
{\bf Proof.}
The proof of this Proposition is implicit within the proof of 2.3 in
\cite {BLS} and can also be obtained as a corollary of Theorem 3.1 below
and Theorem 3.1 in \cite{BLS}. However, for completeness, we outline the
proof
here. In a first step define the map
$$\tilde T \subset (U\times G(d,N))\times \BC^N \buildrel{\tilde F}\over
\longrightarrow \BD_\eta \subset \BC$$
by $\tilde F (x,P,y) = d\tilde f_x (y)$. As $0$ is a general point of  
$f$,
then $\tilde K = \tilde T \cap  \tilde F^{-1}(0)$ is a sub-bundle of
$\tilde
T$ of (complex) codimension 1 and $d \tilde F$ maps the orthogonal
complement
of $\tilde K$ isomorphically over $T(\BD_\eta)$ (see \cite {BLS} for
details).

Now, as we show in the proof of Lemma 3.4 below,
every $\tilde f$ defining an isolated singularity at $0$ in $X$,  
determines a
sub-bundle $Q$ of $T\BC^N \vert_{X-\{ 0 \}}$ everywhere
transversal to ${ker}( d \tilde f)$,
and the restriction of
$d \tilde f$  to $Q$  is an isomorphism between $Q$ and $T(\BD_\eta)$.

This implies that each nowhere-zero vector field  on $\BD_\eta$ lifts
compatibly to a vector field on $X-\{ 0 \}$ and also as a section of
$\tilde T
\vert_X$. The final step is to notice, as we do in Lemma 2.2, that the  
gradient
vector field
$\grad$ can be obtained by lifting such a vector field.
\qed

\medskip

\proclaim{PROPOSITION 2.5} Let $f: (X,0) \to (\BC,0)$  be analytic at
$0$. Then there exists a Zariski open subset $\Omega_f$ of the space of
complex linear forms on $\BC^N$, such that for all $l \in \Omega_f$, the
point $0$ is general for the map
$f+ \lambda l  : X   \to \BC  $, for all sufficiently
small  $\lambda \in \BC^*$.
\endproclaim

This Proposition follows from the proof of
the Morsification Theorem 2.2 in \cite {Le2}. In fact, the deformations  
of
$f$ of the type $f+ \lambda l $, constructed by L\^e in \cite{Le2} have  
$0$ as a
general point and  away from 0
their singularities are quadratic (this latter fact is not needed here).
\qed
\bigskip

\head
{3. THE THEOREM}
\endhead

\proclaim{THEOREM 3.1} Let $ f : (X,0)  \to (\BC,0)$ have an isolated
singularity at $0
\in X$. One has:
$$Eu_X(0)=\left(\sum_i \chi(V_i\cap \Bbb B_\varepsilon\cap
f^{-1}(t_0)) \cdot  Eu_X(V_i) \right) +  Eu_{f,X}(0)\,.$$
\endproclaim

\medskip
Theorem 3.1 is proved using Lemma 3.2 below. To state this Lemma, let us
fix some notations.
We choose $
\epsilon > 0$ sufficiently small  so that  every sphere $\BS_\gamma$ in
$U$
centered at $0$ and radius $ \gamma\leq \epsilon$ intersects  
transversally
every
stratum in $X \setminus \{0 \} $.  For each $ t \in \BC$, set $Y_t :=
f^{-1}(t)$.
Choose   $\eta >0$ small enough  so that for each $t$  in the disk
$\BD_\eta$ of radius $\eta$ around $0 \in \Bbb C$, the hypersurface
$Y_t$ intersects transversally the sphere $\BS_\epsilon$ .
Now choose $\epsilon '$ with $0 < \epsilon' < \epsilon$, and a point
$t_0 \in \BD_\eta$ such that $Y_{t_0}$ does not meet the sphere
$\BS_{\epsilon'}$.
We notice that the strata
$V_i$  intersect $Y_{t_0} :=f^{-1}(t_0)$ transversally and provide a
Whitney
stratification of this space.

\smallskip
\proclaim{LEMMA 3.2} There is a stratified vector field $\,w$ on
$\,X_{\epsilon',\epsilon} = X\cap (\BB_{\epsilon}
\setminus Int(\BB_{\epsilon'}))$ such that:
\item{i)} it coincides with $\grad$ on $\,X\cap  \BS_{\epsilon'}$ and  
its
restriction to $\,X\cap \BS_{\epsilon}$ is  radial;
\item{ii)} it is tangent to $Y_{t_0}$;
\item{iii)} $w$ has only a finite number of zeroes, and they are all
contained in
$Y_{t_0}$;
\item{iv)} at each zero $a$, $w$ is
transversally radial to  the stratum containing $a$ ({\it i.e.}
  it is  transversal to the boundary of a tubular neighbourhood of the
stratum).
\endproclaim

\noindent {\bf Proof of Theorem 3.1} (assuming  Lemma 3.2).
We first
notice that if $\xi$ is a stratified vector field on a neighbourhood of
$\{ 0\}$ in $X$, which is everywhere transversal to a small sphere
$\Bbb S_\varepsilon$, then $\xi$ is homotopic to a radial vector field,
by elementary obstruction theory. Hence to compute the Euler obstruction
using 1.2, it is enough to consider vector fields transversal to
$\Bbb S_\varepsilon$.

The  restriction of the vector field $w$ of 3.2 to
$\partial(X_{\epsilon',\epsilon})$ is a stratified vector field, so
it can be lifted as a section
$\tilde w$ of the Nash bundle $\tilde T$ on
$\nu^{-1}(\partial(X_{\epsilon',\epsilon}))$ by 1.1. Let us denote by
$\text{Obs}(\tilde  w,\nu^{-1}(X_{\epsilon}))$
the obstruction for extending $\tilde w$ to
$\nu^{-1}(X_{\epsilon})$. One has:
$$\text{Obs}(\tilde  w,\nu^{-1}(X_{\epsilon})) =
\text{Obs}(\tilde  w,\nu^{-1}(X_{\epsilon'})) +
\text{Obs}(\tilde  w,\nu^{-1}(X_{\epsilon',\epsilon}))$$

By statement (i) in Lemma 3.2 this formula becomes
$$Eu_X(0) =  Eu_{f,X}(0) + \text{Obs}(\tilde
w,\nu^{-1}(X_{\epsilon',\epsilon})) \,.$$
By statement iii) in the Lemma 3.2, the contribution of  
$\text{Obs}(\tilde
w,\nu^{-1}(X_{\epsilon',\epsilon})) $ is
concentrated on
$\nu^{-1}(Y_{t_0} \cap \BB_{\epsilon})$.
Statements iii) and iv), together with  the
``Theorem of Proportionality'' (\cite {BS}, Th\'eor\`eme 11.1),
  imply that the contribution
of each singularity $x$ of $w$ to $\text{Obs}(\tilde
w,\nu^{-1}(X_{\epsilon',\epsilon})) $  is
$Eu_X(x)$-times
the local Poincar\'e-Hopf index of $w$ at $x$, regarded as a vector  
field
on the
stratum $V_i(x)$. Furthermore (by ii and iv), the sum of the
Poincar\'e-Hopf indices of the restriction of $\,w$ to $V_i \cap  
Y_{t_0}$
is
$\chi(V_i \cap Y_{t_0} \cap \BB_{\epsilon})$.
Theorem 3.1 then follows. \qed

\medskip

The next result follows immediately from Theorem 3.1 and Proposition 2.5
above.

\proclaim{COROLLARY 3.3} Let $f: (X,0) \to (\BC,0)$ have an isolated
singularity at $0$. Choose a  linear form $l\in \Omega_f$ as in 2.5,
and $\lambda \in \BC^*$ so that
$f_{\lambda} :=
f+ \lambda l$ is general at $0$. Let $M_{f,X}$ and $M_{f_{\lambda},X}$  
be
the
Milnor fibres of $f$ and $f_{\lambda}$, respectively, on $X$. Then,
$$Eu_{f,X}(0) = \sum_i \big[\chi (V_i \cap M_{f_{\lambda},X} ) -
\chi (V_i \cap M_{f,X})\big] \cdot Eu_X(V_i) \,$$
\endproclaim

\noindent
{\bf REMARK 3.4} We notice that if  $X \hskip -2pt =\hskip -1pt \BC^d$  
and
$f: (\BC^d,0)  \to
(\BC,0) $ is an analytic map with an isolated singularity at $0$ with  
Milnor
number
$\mu$, then Theorem 3.1 implies,
$$Eu_{f,X}(0) = (-1)^d  \mu \,,$$
because the Euler-Poincar\'e characteristic of the Milnor fibre of $f$
is equal to $1+(-1)^{d-1} \mu$ by \cite {Mi}.

This formula implies the Theorem 7.2 of \cite {Mi}. In fact,
since in this case the variety $X= \BC^d$ is smooth, its Nash transform  
is
$\BC^d$ and the Nash bundle $\tilde T$ is the tangent bundle of $\BC^d$.
Hence, by definition, $Eu_{f,X}(0)$ is the Poincar\'e-Hopf index at $0$
of the gradient vector field
$\overline \nabla f = (\frac {\overline{\partial f}}{{\partial
x_1}},..., \frac {\overline{\partial f}}{\partial x_N}) \,.$ This equals
$(-1)^d$-times the Poincar\'e-Hopf index at $0$ of the vector field
$ \nabla f = (\frac { {\partial f}}{{\partial
x_1}},..., \frac { {\partial f}}{\partial x_N}) \,.$ So we have that
the Milnor number $\mu$ is the degree of $\nabla f$, which is Milnor's
Theorem 7.2.

\bigskip
\head
{4. PROOF OF LEMMA 3.2}
\endhead

The vector field $w$ of Lemma 3.2 is defined as the sum of two
vector
fields that we construct in Lemmas 4.1 and 4.2 below. Let us choose
$\eta'> 0$
sufficiently small with respect to $\eta$, such that the disc
$\BD_{\eta{'}}(t_0)$ of radius $\eta'$, centered in $t_0$, is contained  
in
the
interior of $\BD_{\eta}$ and $f^{-1}(\BD_{\eta{'}}(t_0))$ does not
intersect
the sphere $\BS_{\epsilon'}$.

\proclaim{LEMMA 4.1} There exists a stratified vector field $r$ on
$\Bbb B_{\epsilon} \cap f^{-1}(\BD_\eta)$ satisfying:

\item{i)} The restriction of $r$ to $\BS_{\epsilon} \cap  
f^{-1}(\BD_\eta)$ ,
is tangent to all the fibers $f^{-1}(t)$,
and  is transversal to $\BS_{\epsilon}$, pointing outwards;

\item{ii)} $r$ is tangent to the fiber $f^{-1}(t_0)$, where it has only
isolated
singularities (zeroes). Furthermore, at each zero,
$r$ is transversally radial $($in $f^{-1}(t_0))$ to the stratum
which contains that zero;

\item{iii)}
$r$ is tangent to the fiber $f^{-1}(t)$ for all $t \in
\BD_{\eta{'}}(t_0)$.

\endproclaim

\vglue 10truecm
\centerline{Figure 1.  The vector field $r$}

\medskip
\noindent {\bf Proof.}
  We first
construct $r$ on $f^{-1}(t_0)$ satisfying ii)   and pointing out of
$\BS_{\epsilon}$. This is done using the  technique of M. H. Schwartz,  
as
in the proof of Lemma 3.4 in
\cite {BLS}. The proof is now similar to the proof of Lemma 2.1 in \cite
{BLS}:
by the  fibration theorem of L\^e (see \cite {Le2} Theorem 1.3), the
function
$f$ determines a stratified  locally trivial fibration of
$\Bbb B_{\epsilon} \cap f^{-1}(\BD_\eta)
\setminus f^{-1}(0)$ over
$\BD_{\eta} \setminus  \{0\}$. This induces, in particular, a trivial
fibration over $\BD_{\eta_{'}}(t_0)$. Hence $r$ can be extended,
as a product, to all the fibers over $\BD_{\eta_{'}}(t_0)$, satisfying
statement iii).

Let us choose $\varepsilon''$, with
$\varepsilon'<<\varepsilon''<\varepsilon$ and such that the restriction  
of
$f$ to $f^{-1}(\BD_{\eta})\cap (\BB_\varepsilon\setminus
\text{Int}(\BB_{\varepsilon''}))$ (shaded part of figure 1) is a trivial
fibration.
  Since $\BD_\eta$ retracts to $\BD_{\eta_{'}}(t_0)$, the
vector field $r$ can be extended to
$\BS_{\epsilon} \cap f^{-1}(\BD_\eta)$ being tangent to the fibers of
$f$
and transversal to the sphere  $\BS_{\epsilon}$, pointing outwards

Using a suitable partition of unity, we can extend $r$ as zero in the
complement of a neighbourhood of $\big(f^{-1}(\BD_{\eta})\cap
(\BB_\varepsilon\setminus \text{Int}(\BB_{\varepsilon''}))\big)\cup
f^{-1}(\BD_{\eta_{'}}(t_0))$ (see figure 1).
\hfill
\qed
\bigskip
Notice that statements  i) and ii) in 4.1 imply that for each stratum
$Y_{t_0} \cap V_i \cap B_{\epsilon}$, the sum of the Poincar\'e-Hopf
indices of
the restriction of $r$ is $\chi (Y_{t_0} \cap V_i \cap B_{\epsilon})$.

\proclaim{LEMMA 4.2} There exists a stratified vector field $u$ defined  
on
$(\Bbb B_{\epsilon} \cap f^{-1}(\BD_\eta)) \setminus\{ 0 \}$ and  
satisfying:

\item{i)} $u$ is tangent to $\BS_\varepsilon$;

\item{ii)} its zero set is $f^{-1}(t_0)$, and $u$ is transversally
radial to $f^{-1}(t_0)$ ;

\item{iii)} $u$ is transversal to $X \cap f^{-1}(\partial \BD_\eta)$.

\endproclaim

  \vglue 10truecm
\centerline{Figure 2.  The vector field $u$}
\medskip
\noindent {\bf Proof.} By hypothesis, the restriction of $f$ to every
stratum
$V_i$ (other than $\{ 0 \}$) is regular, hence the kernel of $df$ has
codimension
1 in $V_i$. The Hermitian metric on $\BC^N$ induces a metric on $V_i$  
and
defines a splitting of the tangent bundle $TV_i$ as the sum of the  
tangent
bundle to the fiber and the normal bundle. The derivative
$df$,  restricted to the normal bundle, is an
isomorphism. Therefore, we can lift every vector field on
$\BD_\eta$ as a vector field tangent to $V_i$ and ``orthogonal" to the
fibers $f^{-1}(t)$. Let us denote by $\xi$ the vector
field on $\BD_\eta$, radial from $t_0$, and by
$u_i$ the lifting of $\xi$ to $V_i$
(see Figure 2).

Just as in the definition of the vector field $\grad$,
the Whitney conditions allow us to glue the different $u_i$ in a
   stratified vector field $u$ defined on
   $(\Bbb B_{\epsilon} \cap f^{-1}(\BD_\eta)) \setminus \{ 0
\}$ and satisfying the conditions of the Lemma. \qed

\bigskip
\goodbreak

\noindent {\bf Proof of Lemma 3.2.}

We first define the vector field $w$ on $(\Bbb B_{\epsilon} \cap
f^{-1}(\BD_\eta))
\setminus \text{Int}( \BB_{\epsilon'})$ as the sum (via a partition of
unity) of
the vector fields that we constructed in Lemmas 4.1 and 4.2. We obtain a
stratified vector field satisfying:

\item{i)} it  is
transverse to the
sphere  $ \BS_{\epsilon}$ , pointing outwards ;

\item{ii)} it is tangent to $Y_{t_0}=f^{-1}(t_0)$;

\item{iii)} the singularities of $w$ are all contained in $Y_{t_0}$;

\item{iv)}  on  each stratum   $Y_{t_0} \cap V_i \cap \BB_{\epsilon} $ ,
$w$ only vanishes at a
finite number of points and the sum of the Poincar\'e-Hopf indices of  
the
restriction of $w$ is $\chi (Y_{t_0} \cap V_i \cap \BB_{\epsilon})$;

\item{v)}  at each of its singular points, $w$ is transversally radial
to the stratum that contains this singular point;

\item{vi)} it is transverse to the
boundary  $ f^{-1}(\partial \BD_\eta) $, pointing outwards.

Then we extend
$w$ to all of $\,X_{\epsilon',\epsilon}$ using Theorem 2.3 in \cite  
{BLS}.
To complete the proof of Lemma 3.2, we notice that
restricted to $\BS_{\varepsilon'}$, the vector field
$w$ is homotopic to the gradient vector field $\grad$.
This fact follows from  Lemma 2.2 and the fact that,
on $\BS_{\varepsilon'}$, the vector field $w$ is the vector field
$u$ of Lemma 3.4, which is homotopic to a constant vector field near 0.
\qed
\bigskip
\head
{5. INTERPRETATION IN TERMS OF VANISHING CYCLES}
\endhead

Recall our notation from section 2: $X$ is a $d$-dimensional analytic  
subset of an open subset $U$ of $\Bbb C^N$, $f$ is the restriction to  
$X$ of a complex analytic function $\tilde f:(U, 0)\rightarrow (\Bbb C,  
0)$, and $\{V_i\}$ is a complex analytic Whitney stratification of $X$.

\vskip .2in

In Definition 2.1, we defined the local Euler obstruction of $f$ in  
what we consider to be the most natural way -- as a geometric  
obstruction. We then proved in Theorem 3.1 that, when $f$ has an  
isolated singularity, the local Euler obstruction of $f$ measures the  
defect in $f$ satisfying the local Euler condition.

\vskip .2in

In this section, we use a formal, derived category argument to prove  
that the defect
$$D_{f, X}(0):=Eu_X(0) - \Big(\sum_i\chi\big(V_i\cap B_\epsilon\cap  
f^{-1}(t_0)\big)\cdot Eu_X(V_i)\Big)$$
has a nice interpretation in terms of vanishing cycles, even when the  
singularity is non-isolated.
\vskip .2in

Let $\bold A^\bullet$ be a bounded complex of sheaves of complex vector  
spaces on $X$ which is constructible with respect to $\{V_i\}$.  We  
will need a number of notions from the derived category, including the  
characteristic cycle of $\bold A^\bullet$. The reader is referred in  
[{K-S}].

For each $V_i$, we denote the closure of the conormal variety of $V_i$  
in $U$ by $\overline{T^*_{V_i}U}$. We let $\operatorname{Ch}(\bold  
A^\bullet) =\sum_i m_i\left[\overline{T^*_{V_i}U}\right]$ denote the  
characteristic cycle of $\bold A^\bullet$ in $U$. We let $\psi_f\bold  
A^\bullet$ and $\phi_f\bold A^\bullet$ denote the nearby and vanishing  
cycles of $\bold A^\bullet$ along $f$, respectively.

\vskip .2in

\noindent{\bf LEMMA 5.1}. {\it There exists a complex $\bold A^\bullet$  
on $X$, constructible with respect to $\{V_i\}$, such that  
$\operatorname{Ch}(\bold A^\bullet)  
=\left[\overline{T^*_{X_{\operatorname{reg}}}U}\right]$
}

\vskip .2in

\noindent{\bf Proof}. This is a trivial induction. See, for instance,  
Lemma 3.1 of  [{Ma1}]. \qed

\vskip .2in

\noindent{\bf DEFINITION 5.2}. We call a complex like that of Lemma 5.1  
a {\it characteristic complex for $X$} (constructible with respect to  
$\{V_i\}$).

\vskip .2in

\noindent{\bf THEOREM 5.3}. {\it Let $\bold A^\bullet$ be a  
characteristic complex for $X$.

Then, $D_{f, X}(0)$ equals negative the Euler characteristic of the  
stalk cohomology of $\phi_f\bold A^\bullet$ at the origin, i.e.,
$$
Eu_X(0) = \Big(\sum_i\chi\big(V_i\cap B_\epsilon\cap  
f^{-1}(t_0)\big)\cdot Eu_X(V_i)\Big)-\chi(\phi_f\bold A^\bullet)_0.
$$

}

\vskip .2in

\noindent{\bf Proof}. By the index formula of Brylinski, Dubson, and  
Kashiwara [{BDK}], for any point $x\in X$, $\chi(\bold A^\bullet)_x$ is  
equal to the local Euler obstruction $Eu_X(x)$.

By the distinguished triangle relating the nearby cycles $\psi_f \bold  
A^\bullet$ and the vanishing cycles $\phi_f \bold A^\bullet$, we obtain  
that $\chi(\bold A^\bullet)_0 = \chi(\psi_f \bold A^\bullet)_0 -  
\chi(\phi_f \bold A^\bullet)_0$.

Let $F:=B_\epsilon\cap f^{-1}(t_0)$ denote the Milnor fibre of $f$ at  
the origin. As $F$ transversely intersects the strata $\{V_i\}$, $F$ is  
Whitney stratified by the (finite) collection of strata $\{F\cap  
V_i\}$. The Euler characteristic of the hypercohomology $\chi\Big(\Bbb  
H^*(F;\,\bold A^\bullet_{|_F})\Big)$ is equal to $\sum_i\chi(F\cap  
V_i)\chi(\bold A^\bullet)_{p_i}$, where $p_i$ is a point in the stratum  
$V_i$.

\vskip .1in

Therefore, we obtain  $$\chi(\psi_f \bold A^\bullet)_0 =  
\chi\big(H^*(F; \bold A^\bullet)\big) = \sum_i \chi(V_i\cap B_\epsilon  
\cap f^{-1}(t_0))\cdot\chi(\bold A^\bullet)_{p_i},$$ where $p_i$ is any  
point in $V_i$.

Combining the above steps, we obtain the desired formula.
\qed

\vskip .2in

\noindent{\bf COROLLARY 5.4}. {\it Let $\bold A^\bullet$ be a  
characteristic complex for $X$.

If $f$ has an isolated singularity at $0$, then $(0, d_0\tilde f)$ is  
an isolated point of the intersection  
$\overline{T^*_{X_{\operatorname{reg}}}U}\cap\operatorname{im}d\tilde  
f$, and
$$
Eu_{f, X}(0)=  - \chi(\phi_f \bold A^\bullet)_0 =  
(-1)^d\big(\overline{T^*_{X_{\operatorname{reg}}}U}\ \cdot\  
\operatorname{im}d\tilde f\big)_{(0, d_0\tilde f)},
$$
where this intersection number is equal to the Milnor number of $f$ at  
the origin in the case where $X$ is affine space.
}

\vskip .2in

\noindent{\bf Proof}. The first equality follows from Theorem 3.1. The  
second equality follows from the index formula of  Ginsburg [{Gi}],  
L\^e [{Le1}], and Sabbah [{Sa2}]. That the intersection number equals  
the Milnor number in the affine case is well-known; see, for instance,  
Proposition I.2.19 of [{Ma3}]. \qed

\bigskip 
\head
{6. CALCULATION IN TERMS OF L\^E NUMBERS}
\endhead
 In this section we 
 describe an algebraic method for calculating $\chi(\phi_f\bold  
A^\bullet)_0$, when $\bold A^\bullet$ is a characteristic complex for  
$X$, even when the critical locus of $f$ is non-isolated. In the affine  
case, this method uses L\^e cycles and numbers, as described in [Ma2].  
When $X$ is arbitrary, one must use the more general L\^e-Vogel cycles  
and numbers, as described in [Ma3].

\vskip .1in

In either case, our method of calculation requires a ``generic'' choice  
of coordinates for the ambient, affine space $U$. This choice of  
coordinates is made as follows. Refine the stratification $\{V_i\}$ to  
a stratification $\Cal W:=\{W_j\}$ such that $\Cal W$ satisfies Thom's  
$a_f$ condition and such that $f^{-1}(0)$ is a union of strata of $\Cal  
W$.  While this refinement certainly produces non-isolated $\Cal  
W$-stratified critical points of $f$, this will cause us no difficulty  
as we are considering the non-isolated case. (By the results of [BMM]  
and [P], a Whitney refinement of $\{V_i\}$ such that $f^{-1}(0)$ is a  
union of strata would satisfy Thom's $a_f$ condition; however, we are  
specifically {\bf not} assuming that $\Cal W$ is still a Whitney  
stratification. That one does not require the Whitney conditions on the  
refined stratification makes it easier to produce such refinements.)

Choose the first coordinate $z_1$ so that the hyperplane $z_1^{-1}(0)$  
transversely intersects, in some neighborhood of the origin, all  
positive-dimensional strata of $\{W_j\}$. Then, there is an induced  
stratification (on the germ at the origin) of $X\cap z_1^{-1}(0)$ given  
by $\{W_j\cap z_1^{-1}(0)\}$. We choose $z_2$ so that $z_2^{-1}(0)$  
transversely intersects, in some neighborhood of the origin, all  
positive-dimensional strata of $\{W_j\cap z_1^{-1}(0)\}$. We continue  
in this inductive manner to choose $z_1, \dots, z_{{}_N}$. We call such  
a coordinate choice {\it prepolar} (at the origin).

Prepolar coordinates are not as generic as possible, but they are  
generic enough for our purposes. Being prepolar at the origin implies  
that the coordinates are also prepolar at each point in a neighborhood  
of the origin, and we assume that we are in such a neighborhood  
throughout the remainder of sections 5 and 6. Note that, if we began (before  
we refined $\{V_i\}$) with an $X$ that had a one-dimensional singular  
set at the origin and an $f$ with a one-dimensional stratified critical  
locus, then the choice of prepolar coordinates is trivial.

\vskip .1in

Fix a choice of prepolar coordinates $(z_1, \dots, z_{{}_N})$.

\vskip .2in

\noindent{\bf The affine case}: We will now describe the case where  
$X=U$ is an open subset of affine space. This case will serve as a good  
introduction to the general case.

In this case, we may take $\bold A^\bullet$ to be the constant sheaf  
$\Bbb C_{U}^\bullet$. We write $\Sigma f$ for the ordinary critical  
locus of $f: U\rightarrow \Bbb C$. Our problem is to calculate  
$\chi(\phi_f\Bbb C_{U}^\bullet)_0$.

As described in [Ma2], this Euler characteristic is an alternating sum  
of L\^e numbers $\lambda^k_f$, which are defined as intersection  
numbers of linear subspaces with the L\^e cycles. We need to define  
these devices.

Consider the analytic cycle $[U]$. This cycle can be written as a sum  
of purely $N$-dimensional cycles $$[U]=\Gamma^N_f+\Lambda^N_f,$$ where  
no component of $\Gamma^N_f$ is contained in $\Sigma f$ and every  
component of $\Lambda^N_f$ is contained in $\Sigma f$. Of course, if  
$f$ is not constant on a connected-component of $U$, we have  
$[U]=\Gamma^N_f$.

Now, we define $\Gamma^k_f$ and $\Lambda^k_f$ by downward induction. If  
we have defined the purely $(k+1)$-dimensional cycle $\Gamma^{k+1}_f$,  
then the hypersurface (cycle) $\dsize\Big[V\Big(\frac{\partial  
f}{\partial z_{k+1}}\Big)\Big]$ properly intersects $\Gamma^{k+1}_f$  
inside $U$ (this is a result of the coordinates being prepolar), and  
therefore there is a well-defined, purely $k$-dimensional intersection  
cycle
$$
\Gamma^{k+1}_f \ \cdot\ \Big[V\Big(\frac{\partial f}{\partial  
z_{k+1}}\Big)\Big],
$$
which we can decompose as
$$\Gamma^{k+1}_f \ \cdot\ \Big[V\Big(\frac{\partial f}{\partial  
z_{k+1}}\Big)\Big]\ =:\ \Gamma^k_f+\Lambda^k_f,$$
where no component of $\Gamma^k_f$ is contained in $\Sigma f$ and every  
component of $\Lambda^k_f$ is contained in $\Sigma f$.

The cycle $\Lambda^k_f$ is the {\it $k$-dimensional L\^e cycle of $f$}.  
As our coordinates are prepolar, $\Lambda^k_f$ properly intersects the  
linear subspace $V(z_1, \dots, z_k)$ at the origin, and we define the  
{\it $k$-dimensional L\^e number of $f$ at $0$} to be the intersection  
number
$$
\lambda^k_f(0):= \big(\Lambda^k_f\ \cdot\ V(z_1, \dots, z_k)\big)_0.
$$
Note that $\Lambda^0_f$ is purely $0$-dimensional; thus,  
$\lambda^0_f(0)$ is simply the coefficient of $[0]$ in $\Lambda^0_f$.

\vskip .1in

\proclaim{PROPOSITION 6.1}
If we let $s:={\operatorname{dim}}_0\Sigma f$, the L\^e numbers  
$\lambda^k_f(0)$ are zero if $k>s$, and the Euler characteristic of the  
vanishing cycles is given by
$$
D_{f, X}(0) = -\chi(\phi_f\Bbb C_{U}^\bullet)_0 = \sum_{k=0}^{s}  
(-1)^{N-k}\lambda^k_f(0).
$$
\endproclaim

In the isolated case, $s=0$ and $\lambda^0_f(0)$ is simply the Milnor  
of $f$ at the origin.
Hence, we recover the result stated in Remark 3.4.

\vskip .3in

\noindent{\bf The general case}: We redo the L\^e cycle and number  
construction, but we work in the cotangent bundle to $U$ and with  
conormal varieties.

The cotangent bundle $T^*U@>\pi>>U$ is isomorphic to the trivial bundle  
$U\times\Bbb C^N\rightarrow U$. The choice of the prepolar coordinates  
$(z_1, \dots, z_{{}_N})$ determines a basis $dz_1, \dots, dz_{{}_N}$  
for the cotangent vectors. We use $w_1, \dots, w_N$ for coordinates  
with respect to this basis. We carry out our calculations with  
coordinates $(z_1, \dots, z_{{}_N}, w_1, \dots, w_N)$ on $U\times\Bbb  
C^N$.

Recall that $\tilde f:U\rightarrow\Bbb C$  is our extension of $f$ to  
all of $U$. Note that in our coordinates, the image of $d\tilde f$ is  
given by
$${\operatorname{im}}d\tilde f := V\left(w_1-\frac{\partial\tilde  
f}{\partial z_1}, \dots, w_N-\frac{\partial\tilde f}{\partial  
z_{{}_N}}\right).$$
It will be important to us that the projection $\pi$, restricted to  
${\operatorname{im}}d\tilde f$, induces an isomorphism onto its image;  
the inverse map sends a point $p$ to $\dsize\left(p,{  
\frac{\partial\tilde f}{\partial z_1}}_{|_p}, \dots, {  
\frac{\partial\tilde f}{\partial z_{{}_N}}}_{|_p} \right)$.

\vskip .2in

The cycle $\Big[\overline{T^*_{X_{\operatorname{reg}}}U}\Big]$ can be  
written as a sum of purely $N$-dimensional cycles  
$$\Big[\overline{T^*_{X_{\operatorname{reg}}}U}\Big]=\widehat\Gamma^N_f+ 
\widehat\Lambda^N_f,$$ where no component of $\widehat\Gamma^N_f$ is  
contained in ${\operatorname{im}}d\tilde f$ and every component of  
$\widehat\Lambda^N_f$ is contained in  ${\operatorname{im}}d\tilde f$.

Now, we define $\widehat\Gamma^k_f$ and $\widehat\Lambda^k_f$ by  
downward induction. If we have defined the purely $(k+1)$-dimensional  
cycle $\widehat\Gamma^{k+1}_f$, then the hypersurface  
$\dsize\Big[V\Big(w_{k+1}-\frac{\partial \tilde f}{\partial  
z_{k+1}}\Big)\Big]$ properly intersects $\widehat\Gamma^{k+1}_f$ inside  
$U$, and therefore there is a well-defined, purely $k$-dimensional  
intersection cycle
$$
\widehat\Gamma^{k+1}_f \ \cdot\ \Big[V\Big(w_{k+1}-\frac{\partial  
\tilde f}{\partial z_{k+1}}\Big)\Big],
$$
which we can decompose as
$$\widehat\Gamma^{k+1}_f \ \cdot\ \Big[V\Big(w_{k+1}-\frac{\partial \tilde f}{\partial  
z_{k+1}}\Big)\Big]\ =:\ \widehat\Gamma^k_f+\widehat\Lambda^k_f,$$
where no component of $\widehat\Gamma^k_f$ is contained in  
${\operatorname{im}}d\tilde f$ and every component of  
$\widehat\Lambda^k_f$ is contained in ${\operatorname{im}}d\tilde f$.

As the $\widehat\Lambda^k_f$ are contained in  
${\operatorname{im}}d\tilde f$, the projection, $\pi$, maps each  
$\widehat\Lambda^k_f$ isomorphically onto a cycle in $U$; we let  
$\Lambda^k_f:=\pi\big(\widehat\Lambda^k_f\big)$ (this is the proper  
projection of a cycle, and is frequently denoted by $\pi_*$). We refer  
to $\Lambda^k_f$ as the {\it $k$-dimensional L\^e-Vogel cycle}.

Note that in the affine case,  
$\Big[\overline{T^*_{X_{\operatorname{reg}}}U}\Big]=U\times \{0\}$, and  
the L\^e-Vogel cycles coincide with the L\^e cycles.

\vskip .2in

Now, exactly as before, $\Lambda^k_f$ properly intersects the linear  
subspace $V(z_1, \dots, z_k)$ at the origin, and we define the {\it  
$k$-dimensional L\^e-Vogel number of $f$ at $0$} to be the intersection  
number
$$
\lambda^k_f(0):= \big(\Lambda^k_f\ \cdot\ V(z_1, \dots, z_k)\big)_0.
$$

\vskip .2in
\proclaim{THEOREM 6.2}
If we let  
$s:={\operatorname{dim}}_0\,\pi\Big(\overline{T^*_{X_{\operatorname{reg} 
}}U}\ \cap\ {\operatorname{im}}d\tilde f\Big)$, the L\^e-Vogel numbers  
$\lambda^k_f(0)$ are zero if $k>s$, and the Euler characteristic of the  
vanishing cycles is given by
$$
D_{f, X}(0) = -\chi(\phi_f\bold A^\bullet)_0 = \sum_{k=0}^{s}  
(-1)^{d-k}\lambda^k_f(0).
$$
\endproclaim
\vskip .2in

Note that when $s=0$, the only L\^e-Vogel number which is possibly  
non-zero is $\lambda^0_f(0)$, and $\lambda^0_f(0)=  
\big(\overline{T^*_{X_{\operatorname{reg}}}U}\ \cdot\  
\operatorname{im}d\tilde f\big)_{(0, d_0\tilde f)}$. Therefore, we  
recover the result of Corollary 5.4.


\bigskip


\Refs

\widestnumber\key{BMM}

\ref
\key Br  \by J.P. Brasselet
\paper Existence des classes de Chern en th\'eorie bivariante
\jour Ast\'erisque,
\vol 101-102 \yr 1983 \pages 7--22
\endref


\ref
\key BS  \by J.P. Brasselet et M.H. Schwartz
\paper Sur les classes de Chern d'un ensemble analytique complexe
\jour Ast\'erisque  \vol 82-83
\yr 1981  \pages 93-147
\endref

\ref
\key BLS  \by J.P. Brasselet,  L\^e D. T., J. Seade
\paper Euler obstruction and indices of vector fields
\jour Topology   \vol 39  \yr 2000   \pages 1193-1208
\endref


\ref \key{BMM}  \by J. Brian\c con, P. Maisonobe, and M. Merle   \paper  
Localisation de
syst\`emes diff\'erentiels, stratifications de Whitney et condition de  
Thom   \yr 1994
\jour Invent. Math.\vol 117 \pages 531--550 \endref

\ref \key{BDK}  \by J. Brylinski, A. Dubson, and M. Kashiwara    \paper  
Formule de l'indice
pour les modules holo\-nomes et obstruction d'Euler locale   \yr 1981   
\jour C.R. Acad. Sci.,
S\'erie A  \vol 293 \pages 573--576
\endref


\ref
\key FM \by W. Fulton and R. MacPherson
\book Categorical Framework for the study of singular Spaces
\publ Memoirs of Amer. Math. Soc.
\vol 243 \yr 1981
\endref

\ref \key{Gi}  \by V. Ginsburg    \paper Characteristic Varieties and  
Vanishing Cycles    \jour
Invent. Math. \vol 84 \yr 1986 \pages 327--403  \endref  


\ref \key GSV  \by  X. G\'omez-Mont, J. Seade and A. Verjovsky
\paper The index of a holomorphic flow with an isolated singularity
\jour  Math. Ann
\vol  291  \yr 1991  \pages  737-751
\endref


\ref
\key GM \by M. Goresky and R. MacPherson
\book Stratified Morse Theory
\publ
Springer Verlag \yr 1987
\endref

\ref \key Go \by V. Goryunov \paper
Simple functions on space curves
\jour Funct.  Anal.  Appl. \vol 34, no. 2  \yr  2000 \pages 129--132
\endref

\ref \key{K-S}   \by M. Kashiwara and P. Schapira     \book Sheaves on  
Manifolds \yr 1990
\publ Grund. math. Wiss. 292, Springer - Verlag   \endref 

\ref
\key Le1  \by  L\^e D. T.
\paper Le concept de singularit\'e isol\'ee de fonction analytique
\jour Adv. Stud. Pure Math. \vol 8
\yr 1986 \pages 215-227
\publ North Holland
\endref


\ref
\key Le2  \bysame \paper  Complex analytic functions with
isolated singularities \jour J. Algebraic Geometry \vol 1
\yr 1992 \pages 83-100
\endref

\ref
\key  MP   \by R. MacPherson \paper Chern classes for singular varieties
\jour Ann. of Math \vol 100  \yr
1974 \pages 423-432
\endref

\ref \key{Ma1}   \by D. Massey   \paper Hypercohomology of Milnor  
Fibres \jour Topology
\vol 35\pages 969--1003\yr 1996    \endref 

\ref \key{Ma2}  \bysame     \book L\^e Cycles and Hypersurface  
Singularities
\yr 1995 \publ Springer-Verlag \bookinfo Lecture Notes in Mathematics,  
vol. 1615
    \endref


\ref
\key Ma3  \bysame
\book Numerical Control over Complex Analytic Singularities
\toappear
\bookinfo Memoirs of the AMS
\publ AMS
\endref

\ref \key MS \by D. Mond and D. Van Straten  \paper
Milnor number equals Tjurina number for functions on space curves \jour
J. London Math. Soc. \vol 63 \yr 2001 \pages 177--187
\endref

\ref
\key Mi  \by J. Milnor  \book Singular Points of Complex Hypersurfaces
\publ Annals
of Mathematics Studies {\bf 61}, Princeton University Press, Princeton  
\yr
1968
\endref

\ref \key{P}  \by A. Parusi\'nski      \pages 99--108 \paper Limits of  
Tangent Spaces to Fibres and the $w_f$ Condition
\yr 1993  \jour Duke Math. J. \vol 72    \endref

\ref
\key  Sa1   \by C. Sabbah
\book Espaces conormaux bivariants
\publ Th\`ese, Universit\'e de Paris VII   \yr
1986
\endref

\ref \key{Sa2}  \bysame      \pages 161--192 \paper Quelques remarques  
sur la g\'eom\'etrie
des espaces conormaux  \yr 1985
\jour Ast\'erisque  \vol 130 \endref \vskip .1in

\ref
\key  Sc   \by M.-H. Schwartz  \book Champs radiaux sur une  
stratification
analytique
complexe  \publ  Tra\-vaux en cours {\bf 39}, Hermann \yr 1991
\endref

\ref
\key Scu   \by J. Sch\"urmann   \paper A short proof of a formula of
Brasselet,
L\^e and Seade
\paperinfo Preprint 2001, math. AG/ 0201316
\yr    \pages
\endref


\ref
\key Wh  \by H. Whitney   \paper  Local properties of Analytic Varieties
\jour Symposium in honor of M. Morse, Princeton Univ. Press, edited by  
S.
Cairns
\yr 1965
\endref

\endRefs
\enddocument
\end